\documentclass[11pt]{amsart}
\usepackage{mathrsfs}
\usepackage{amsfonts}
\usepackage{epsfig}
\usepackage{amsmath,amssymb,amsthm}

\begin{document}

\newtheorem{thm}{Theorem}
\newtheorem{prop}[thm]{Proposition}
\newtheorem{lem}[thm]{Lemma}
\newtheorem{cor}[thm]{Corollary}
\newtheorem{conj}[thm]{Open Problem}

\theoremstyle{definition}

\newtheorem{defn}[thm]{Definition}
\newtheorem{notation}[thm]{Notation}
\newtheorem{example}[thm]{Example}


\title[A uniqueness result on ordinary differential equations]
{A uniqueness result on ordinary differential equations with
  singular coefficients}
\author[]{Yifei Pan}
\address{Yifei Pan\\
Department of Mathematical Sciences \\ Indiana University -
Purdue University Fort Wayne \\ Fort Wayne, IN 46805-1499}
\address
{School of Mathematics and Informatics \\ 
Jiangxi Normal University, Nanchang, China}
\email{pan@ipfw.edu}

\author[]{Mei Wang}
\address{Mei Wang\\
Department of Statistics\\ University of Chicago\\
Chicago, IL 60637}
\email{meiwang@galton.uchicago.edu}

%
\begin{abstract}
  We consider the uniqueness of solutions of ordinary differential
  equations where the coefficients may have singularities.  We derive
  upper bounds on the the order of singularities of the coefficients
  and provide examples to illustrate the results.
\end{abstract}
\subjclass[2000]{34A12, 65L05}
\maketitle
%
%
%
\section{Results and examples}
%
%
Classical results on the existence and uniqueness of ordinary
differential equations are mostly concerned with continuous
coefficients (ref. \cite{TenenbaumPollard}).  Here we consider the
uniqueness of ordinary differential equation solutions of coefficients
with singularities.  We study upper bounds on the the order of
singularities of the coefficients that guarantee the uniqueness of the
solution.

\vspace{3mm} Main theorems are stated below.  Two examples are given
to illustrate and to address the sharpness aspect of the results.
Proofs are provided in the subsequent section.

\vspace{3mm}
%
%
%
\begin{thm} \label{UC2-thm-e-bound} Let $f(x)\in\mathcal
  C^\infty(-a,a)$ be a solution (real or complex) of
\begin{equation}  \label{UC2-ODE}
y^{(n)} + a_{n-1}(x,y)y^{(n-1)} + \cdots + a_0(x,y) y = 0,
\qquad x\in(-a,a), ~a>0
\end{equation}
with initial conditions
$$f(0)=f'(0)=\cdots = f^{(n-1)}(0)=0. $$
If
\begin{equation} \label{UC2-e-bound}
\lim_{x\to 0}|x|^{n-k}\left|a_k(x,y)\right| 
\le \frac{1}{e}, \qquad k=0,1,\cdots,n-1,
\end{equation}
where $e$ is the Euler's number, then there exists $\delta>0$ such
that
$$ f\equiv 0 \quad on\quad [-\delta,\delta].$$
\end{thm}

\vspace{3mm}\noindent{\bf Remark:}  For fixed $n$, the inequality
(\ref{UC2-e-bound}) can be relaxed to
\begin{equation} \label{UC2-Bn-bound}
\lim_{x\to 0}|x|^{n-k}\left|a_k(x,y)\right|
< \frac{1}{B_n}, \qquad k=0,1,\cdots,n-1, 
\quad B_n=\sum_{k=0}^{n-1}\frac{1}{k!}.
\end{equation}
%

%
%
%
%
\vspace{5mm}
\begin{cor} \label{UC2-cor-nearzero} Let $f(x)\in\mathcal
  C^{\infty}(-a,a)$ be a solution of (\ref{UC2-ODE}) with initial
  conditions
$$ f(0)=f'(0)=\cdots = f^{(n-1)}(0)=0. $$
If
$$ |a_k(x,y)|=o\left(\frac{1}{|x|^{n-k}}\right)
\qquad as \quad x\to 0,\qquad k=0,1,\cdots, n-1,$$
then there exists $\delta>0$ such that 
$$ f\equiv 0 \quad on\quad [-\delta,\delta].$$
\end{cor}

%
%
%
%

\vspace{3mm}
\begin{cor} \label{UC2-cor-bddcoeff} Let $f(x)\in\mathcal
  C^{\infty}(-a,a)$ be a solution of (\ref{UC2-ODE}) with initial
  conditions
$$ f(0)=f'(0)=\cdots = f^{(n-1)}(0)=0. $$
If
$$ |a_k(x,y)|\le M
\qquad as \quad x\to 0,\qquad k=0,1,\cdots, n-1$$
for some $M>0$, then exists $\delta>0$ such that
$$ f\equiv 0 \quad on\quad [-\delta,\delta].$$
\end{cor}

\vspace{3mm}
%
%
%
\begin{example} \label{UC2-ex-notC2} ~This example shows that the
  uniqueness in Theorem \ref{UC2-thm-e-bound} may not be true for
  solutions not sufficiently smooth.  For $\alpha\in (0,1)$, the
  function
\begin{eqnarray*}
y =
\begin{cases}
\quad x^{\alpha} \sin (x),& x\in[0,\infty)\\
(-x)^{\alpha} \sin (-x),& x\in(-\infty,0)
\end{cases}
\end{eqnarray*}
satisfies the differential equation 
\begin{equation} \label{madeup-eq}
y'' - \frac{2\alpha}{x} y' + 
\left(1+\frac{\alpha^2 + \alpha}{x^2}\right) y =0
\qquad \text{\it with}\qquad  y(0)=y'(0)=0.
\end{equation}
Let $\alpha=\frac{1}{2e}$.  Then condition (\ref{UC2-e-bound}) in
Theorem \ref{UC2-thm-e-bound} is satisfied (for $n=2$):
$$\lim_{x\to 0}|x|\left|a_1(x,y)\right|=\frac{1}{e}, 
\qquad
\lim_{x\to
  0}|x|^2\left|a_0(x,y)\right|=\frac{1}{2e}\left(1+\frac{1}{2e}\right) 
 <\frac{1}{e}.
$$
But $y\not\equiv 0$.  Thus solutions to equation (\ref{madeup-eq}) are
not unique. Notice that $y\in\mathcal C^{1,\alpha}$ (first derivative
of H\"older continuity of order $\alpha$), $y \not\in \mathcal
C^\infty$.  The example also shows that the non-uniqueness cannot be
remedied by using a smaller bound in (\ref{UC2-e-bound}), because for
any given $\varepsilon>0$, we may choose $\alpha<\varepsilon/2$ such
that
\begin{equation*} 
\lim_{x\to 0}|x|^{2-k}\left|a_k(x,y)\right|
\le\max_\alpha\{2\alpha, \alpha^2+\alpha\}<\varepsilon, \qquad k=0,1.
\end{equation*}
\end{example}

\vspace{3mm}
%
%
%
%
\begin{example} \label{UC2-ex-Bessel2} This example shows that a bound
  in condition (\ref{UC2-e-bound}) in Theorem \ref{UC2-thm-e-bound} is
  necessary.  Consider the Bessel differential equation
  (ref. \cite{watson})
\begin{equation}\label{Bessel-eq}
y'' + \frac{1}{x} y' + \left(1-\frac{\nu^2}{x^2}\right) y =0.
\end{equation}
A real solution can be of the form 
$$y_\nu(x)=\sum_{k=0}^\infty
\frac{(-1)^k}{k!~\Gamma(k+\nu+1)}
\left(\frac{x}{2}\right)^{2 k + \nu}
=x^{\nu}g(x) $$
where $g(x)$ is real analytic, $g(0)\not=0$.  Let $\nu=m\ge 2$ be an
integer.  Then
$$y_m(x)=x^m g(x)\in\mathcal C^\infty$$ 
is a solution to (\ref{Bessel-eq}) with $y_m'(x)=m x^{m-1}g(x)+x^m
g'(x)$ and $y_m(0)=y_m'(0)=0$. But $y_m(x)\not\equiv 0$.  Thus
solutions to equation (\ref{Bessel-eq}) are not unique.  Notice that
the only assumption not satisfied in Theorem \ref{UC2-thm-e-bound} is
Condition (\ref{UC2-e-bound}):
\begin{eqnarray*}
\lim_{x\to 0}|x|^n\left|a_0(x,y)\right|
&=&
\lim_{x\to 0} |x|^2\left|1-\frac{m^2}{x^2}\right|=m^2
>\frac{1}{e}. 
\end{eqnarray*}

\end{example}

\vspace{3mm}
%
%
%
%
\begin{example}
In the case of Cauchy-Euler or
equi-dimensional equations,
\begin{equation}\label{UC2-ex-CauchyEuler}
x^n y^{(n)} + a_{n-1}x^{n-1}y^{(n-1)} + \cdots + a_0 x y = 0,
\qquad x\in(-a,a)
\end{equation}
where $a_k$'s are constants, Condition (\ref{UC2-e-bound}) is
simplified to 
$$|a_k|<\frac{1}{e}, \qquad k=0,\cdots, n-1.$$
For $n=2$, the solutions for (\ref{UC2-ex-CauchyEuler}) have the forms
$y = c_1 x^\alpha + c_2 x^\beta,~y=c_1 x^\alpha\ln(x)+ c_2 x^\beta$ or
$y = c_1 x^\alpha\cos(\beta\ln(x)) + c_2
x^\alpha\sin(\beta\ln(x))$. These solutions do not fall into the
categories described in Example \ref{UC2-ex-notC2} or Example
\ref{UC2-ex-Bessel2}.
\end{example}

\vspace{5mm}
%
\section{Proofs}
%
%
We need our previous result (\cite{flatpaper}, Theorem 5)
which is stated here as a lemma.
\begin{lem}\label{UC1-thm5-UC2lem} Assume $f$ (real or complex) is in
  $\mathcal C^\infty(a,b), ~0\in(a,b)$, and for $n\ge 2$ and some
  constant $C$,
\begin{equation} \label{UC2-fnbound-reverse} 
|f^{(n)}(x)|\le C\sum_{k=0}^{n-1}\frac{|f^{(k)}(x)|}{|x|^{n-k}},
\qquad x\in(a,b).
\end{equation}
Then
$$f^{(k)}(0)=0, \quad \forall k\ge 0\qquad {\rm implies}\qquad
 f\equiv 0.$$ 
\end{lem}

\vspace{5mm} First we prove a lemma that provides an upper bound on
the vanishing order of $f$ near $0$ when $f\not\equiv 0$.

%
%
%
%
\vspace{3mm}\begin{lem} \label{UC2-lem-f0-reverse} Assume $f(x)\in
  \mathcal C^\infty(a,b)~,0\in(a,b)$, and (\ref{UC2-fnbound-reverse})
  holds for $n\ge 2$ and some constant $C$.  If $f\not\equiv 0$ on
  $(a,b)$, then at $x=0$, $f$ is of finite vanishing order $N$,  
$$ N \le B_n C + n -1,\qquad B_n=\sum_{k=0}^{n-1}\frac{1}{k!}, $$
i.e., there exists $N>0$ such that for $x$ near $0$,
$$f(x)=a_N x^N + O\left(x^{N+1}\right).$$
\end{lem}

\vspace{5mm}

%
%

\begin{proof} 
%
  When $f\not\equiv 0$, by Lemma \ref{UC1-thm5-UC2lem}, there must
  exist $N$ such that
$$f^{(j)}(0)=0, \quad \forall j<N,~j\ge 0,\qquad {\rm and}\qquad
 f^{(N)}(0)=a_N\not=0.$$ 
Since $f(x)\in \mathcal C^\infty(a,b)$, Taylor's theorem yields
$$f(x)=a_N x^N + O\left(x^{N+1}\right).  $$
If $N\ge n-1$, then 
\begin{eqnarray*}
\frac{|f^{(k)}(x)|}{|x|^{n-k}}
&=&\frac{\left|N(N-1)\cdots(N-k+1) a_N x^{N-k} + O(x^{N-k+1})\right|}
{|x|^{n-k}}\\
&=&N(N-1)\cdots(N-k+1) |a_N x^{N-n}| + O\left(|x|^{N-n+1}\right)
\end{eqnarray*}
for $k=1,2,\cdots,n-1$.  By (\ref{UC2-fnbound-reverse}), for $x\in(a,b)$,
as $x$ approach $0$,
\begin{eqnarray*}
|f^{(n)}(x)|\quad
&=&N(N-1)\cdots(N-n+1) |a_N x^{N-n}| + O(|x|^{N-n+1})\\
\le~ C\sum_{k=0}^{n-1}\frac{|f^{(k)}(x)|}{|x|^{n-k}}
&=& C\left(1+\sum_{k=1}^{n-1}  
N(N-1)\cdots(N-k+1)\right) \left|a_N x^{N-n}\right| 
+ O\left(|x|^{N-n+1}\right).
\end{eqnarray*}
If $N\ge n-1$, dividing both sides by
$N(N-1)\cdots(N-n+2)|a_N x^{N-n}|$ we obtain
\begin{eqnarray*}
N-n+1 + O(|x|)
\le C~\frac{1+\sum_{k=1}^{n-1}
N(N-1)\cdots(N-k+1)}{N(N-1)\cdots(N-n+2)} + O(|x|).
\end{eqnarray*}
Letting $x\to 0$,
\begin{eqnarray*}
N\!-\!n+1 
&\le& 
C~\frac{1+\sum_{k=1}^{n-1}
N(N-1)\cdots(N-k+1)}{N(N-1)\cdots(N-n+2)} \\
&=&
C~\frac{1+N + N(N-1)+ \cdots+ N(N-1)\cdots(N-n+2)}
{N(N-1)\cdots(N-n+2)}\\
&=&\!\!
C\!\left(\!\!\frac{1}{N(N\!-\!1)\cdots(N\!-\!n+2)} 
+\frac{1}{(N\!-\!1)\cdots(N\!-\!n+2)} + \cdots
+\frac{1}{N\!-\!n+2} +1\!\!\right)\\
&\le & C\left(\frac{1}{(n-1)!}+\frac{1}{(n-2)!}+\cdots 
+\frac{1}{2!} + \frac{1}{1!}+1\right) =  C B_n.
\end{eqnarray*}
Notice that the last inequality achieves equality when $N=n-1$.
Thus when $N\ge n-1$, the order of $f(x)=a_N x^N + O(x^{N+1})$ satisfies 
$n-1\le N\le B_n C+ n-1$.
Combining with the case of $N < n-1$, we obtain
$$N~ \le ~B_n C ~+~ n-1.$$
This completes the proof of Lemma \ref{UC2-lem-f0-reverse}.
\end{proof}

\vspace{3mm} Next, we consider a proposition slightly more general
than Corollary \ref{UC2-cor-nearzero}.
%
%
%
%

\vspace{3mm}
\begin{prop} \label{UC2-prop1-finite} Let $f\in\mathcal
  C^{\infty}(-a,a)$ be a solution of (\ref{UC2-ODE}) such that
\begin{equation} \label{UC2-prop1-coeffbound} 
|a_k(x,y)|=O\left(\frac{1}{|x|^{n-k}}\right)
\qquad as \quad x\to 0,\qquad k=0,1,\cdots, n-1.
\end{equation}
If
\begin{equation}  \label{UC2-prop1-diffbound}
f^{(k)}(0)=0,\quad \forall k\le B_n C_n + n-1,
\end{equation}
where
$$ C_n = \limsup_{0\le k \le n-1,~x\to 0}
\left\{|a_k(x,y)| |x|^{n-k}\right\}, \qquad
B_n=\sum_{k=0}^{n-1}\frac{1}{k!},$$
then there exists $\delta>0$ such that $f\equiv 0$ on
$[-\delta,\delta]$.
\end{prop}

\vspace{1mm}
%
%
%
\begin{proof} 
It follows from the differential equation (\ref{UC2-ODE}) that 
$$\left|f^{(n)}(x)\right| 
\le \sum_{k=0}^{n-1}|a_k(x,y)| \left|f^{(k)}(x)\right|, \qquad
\forall x \in (-a,a).$$  
Then
$$  C_n=\max_{0\le k \le n-1}c_k, \qquad\text{with}\quad
c_k=\limsup_{x\to 0} \left\{|x|^{n-k} |a_k(x,y)|\right\},
$$
and $c_k$'s are finite by Assumption (\ref{UC2-prop1-coeffbound}).
Therefore, for any given $\varepsilon>0$, there exists $\delta>0$ such
that
$$\left|f^{(n)}(x)\right| 
\le (C_n +\varepsilon)
\sum_{k=0}^{n-1}\frac{\left|f^{(k)}(x)\right|}{|x|^{n-k}}, \qquad
\forall x \in [-\delta, \delta]. $$
If $f\not\equiv 0$ on $[-\delta,\delta]$, we would have
$f^{(N)}(0)\not=0$ for some $N\le B_n(C_n+\varepsilon)+n-1$ by Lemma
\ref{UC2-lem-f0-reverse}, and the arbitrariness of $\varepsilon$
would imply $f^{(N)}(0)\not=0$ for some $N\le M$, where $M=\lfloor B_n
C_n + n - 1\rfloor$ is the largest integer $\le B_n C_n + n - 1$.
However Condition (\ref{UC2-prop1-diffbound}) implies
$f^{(k)}(0)=0,~\forall k\le M$. Hence we must have $f\equiv 0$ on
$[-\delta,\delta]$ for some $\delta>0$.  This completes the proof of
Proposition \ref{UC2-prop1-finite}.
%
\end{proof}  
%

\vspace{3mm}\noindent{\bf Remark:}
Notice that Example \ref{UC2-ex-Bessel2} satisfies Condition
(\ref{UC2-prop1-coeffbound}) in Proposition \ref{UC2-prop1-finite}:
\begin{equation} \label{UC2-ex-Besselbound}
|a_0(x,y)|
=
\left|1-\frac{m^2}{x^2}\right|=O\left(\frac{1}{|x|^{n-0}}\right),\quad
|a_1(x,y)|
=
\left|\frac{1}{x}\right|=O\left(\frac{1}{|x|^{n-1}}\right)
\end{equation}
as $x\to 0$ for $k=0,1$ ($n=2$).  However the uniqueness does not hold
because Condition (\ref{UC2-prop1-diffbound}) is not satisfied:
$y_m^{(m)}\not=0$, where $m<M= B_n C_n + 1, ~C_n=m^2$.

\vspace{5mm} The proof of Theorem \ref{UC2-thm-e-bound} follows from
Proposition \ref{UC2-prop1-finite}, as stated below.
%
%
%
%

\begin{proof} By the assumption in Theorem \ref{UC2-thm-e-bound},
  $C_n=\frac{1}{e}$. Since
$$ B_n C_n + n-1 =B_n\frac{1}{e} + n-1 < e~\frac{1}{e} + n-1=n,$$
the initial conditions $f^{(k)}(0)=0,~\forall k<n$ imply
$$f^{(k)}(0)=0, \qquad\forall k\le B_n C_n + n-1. $$
Therefore $f\equiv 0$ on $|x|\le
\delta$ for some $\delta>0$ by the result in Proposition
\ref{UC2-prop1-finite}.  This completes the proof of Theorem
\ref{UC2-thm-e-bound}.
\end{proof}

%
%
%
%
\vspace{5mm} Similarly, Corollary \ref{UC2-cor-nearzero} follows immediately.
\begin{proof} By the assumption, $C_n=0, ~B_n C_n+n-1 = n-1$.  Since
  $f^{(k)}(0)=0, ~\forall k\le n-1$, the result of Corollary
  \ref{UC2-cor-nearzero} follows from Proposition
  \ref{UC2-prop1-finite}.
\end{proof}

\vspace{9mm}

\vspace{5mm}


\begin{thebibliography}{[10]}

\bibitem{TenenbaumPollard} {\sc Tenenbaum, M} and {\sc Pollard, H.} (1963)
  {\sl Ordinary Differential Equations}, Harper \& Row. New York.

\bibitem{flatpaper} {\sc Y. Pan} and {\sc M. Wang} (2008) {\sl When is a
  function not flat?}  J. Math. Anal. App. {\bf 340}(1), 536-542.

\bibitem{watson} {\sc G. N. Watson} (1944) {\sl A Treatise on the
    Theory of Bessel Functions}, 2nd ed.  Cambridge University press.  

\end{thebibliography}
\end{document}